\documentclass[12pt]{article}
\usepackage[T1]{fontenc}
\usepackage[utf8]{inputenc}
\usepackage{lmodern}
\usepackage{amsmath,amsthm,amsfonts,amssymb}
\usepackage{fullpage}
\usepackage[shortlabels]{enumitem}
\usepackage{microtype}
\usepackage[hidelinks,hyperfootnotes=false]{hyperref}
\hypersetup{pdftitle={A proof of Minh's conjecture}}

\newtheorem{thm}{Theorem}[section]
\newtheorem{lem}[thm]{Lemma}
\newtheorem{prop}[thm]{Proposition}
\theoremstyle{definition}
\newtheorem{ex}[thm]{Example}

\numberwithin{equation}{section}
\makeatletter
\let\c@figure\c@equation
\let\c@table\c@equation

\makeatother
\allowdisplaybreaks[1]

\setlist{leftmargin=3\parindent,labelindent=3\parindent}
\setlist[enumerate]{
  leftmargin=3\parindent,
  align=left,
  labelwidth=3\parindent,
  labelsep=0pt
}
\setlist[enumerate,1]{label={\normalfont (\arabic*)},ref={\arabic*}}

\DeclareMathOperator{\reg}{reg}
\DeclareMathOperator{\odd}{odd}
\newcommand{\kk}{\Bbbk}
\newcommand{\NN}{\mathbb{N}}
\newcommand{\sfq}{\tau_{\mathrm f}}
\newcommand{\arxiv}[1]{\href{https://arxiv.org/abs/#1}{\emph{arXiv:#1}}}

\title{A proof of Minh's conjecture}
\author{
Hao Chen
~~~~~~~
Yulin Hu
}
\date{} 

\begin{document}
\maketitle

\begin{abstract}
We prove Minh’s conjecture that, over any field, the ordinary and symbolic powers of the edge ideal of a finite simple graph have the same Castelnuovo–Mumford regularity at every positive exponent. We also prove that each ordinary power has the same regularity as its integral closure.
\end{abstract}

\begingroup
\renewcommand{\thefootnote}{}
\footnotetext{School of Mathematical Sciences, Soochow University,
Suzhou 215006, China.}
\endgroup

\section{Introduction}

Let $G=(V,E)$ be a finite simple graph and let $S=\kk[x_v:v\in V]$ be the standard graded polynomial ring over a field $\kk$. The edge ideal of $G$ is $I(G)=(x_ux_v:uv\in E)\subseteq S$. For an integer $s\ge1$, write $I(G)^s$ for its ordinary power, $I(G)^{(s)}$ for its symbolic power, and $\overline{I(G)^s}$ for the integral closure of its ordinary power. Simis, Vasconcelos and Villarreal proved that $I(G)^s=I(G)^{(s)}$ for every $s$ precisely when $G$ is bipartite \cite[Theorem~5.9]{SVV}. If $C$ is an odd cycle of length $2r+1$ in $G$, its vertex product $\prod_{v\in V(C)}x_v$ lies in $I(G)^{(r+1)}$, since every vertex cover meets the cycle in at least $r+1$ vertices, but its degree is too small for $I(G)^{r+1}$.

Minh conjectured that the ordinary and symbolic powers nevertheless have the same Castelnuovo--Mumford regularity, denoted by $\reg$, at each exponent \cite[Conjecture~A]{MNP}. A.~Kumar and R.~Kumar asked whether the integral closure of the ordinary power also has this regularity \cite[Question~1, first part]{KK}.

\begin{thm}\label{thm:main}
For every finite simple graph $G$, every field $\kk$, and every integer $s\ge1$, the edge ideal $I(G)\subseteq\kk[x_v:v\in V(G)]$ satisfies
\[
 \reg\bigl(I(G)^s\bigr)
 =\reg\bigl(\overline{I(G)^s}\bigr)
 =\reg\bigl(I(G)^{(s)}\bigr).
\]
\end{thm}

For the rest of the paper, we write $I=I(G)$ for the edge ideal of the fixed graph $G$.

For a graph of odd girth $2r+1$ (the length of its shortest odd cycle), Kumar and Kumar proved Minh's equality at $s=r+1$, the first exponent at which the ideals differ \cite[Theorems~4.3 and~4.5]{KK}. For arbitrary graphs, Minh, Nam, Phong, Thuy and Vu proved the equality for $s=2,3$ \cite[Theorem~1.1]{MNP}. Their comparison uses the radicals of colon ideals at the multidegrees that realize regularity. Extending it requires control of both the additional generators of the symbolic power and the resulting colon ideals as $s$ increases. Minh and Vu proved the ordinary-power and integral-closure comparison for $s\le4$ \cite[Theorem~1.1]{MV}.

For all exponents, Gu, Hà, O'Rourke and Skelton treated odd cycles \cite[Theorem~5.3]{Gu}, Seyed Fakhari treated unicyclic graphs \cite[Theorem~3.9]{FakUni}, and Muta treated simplicial graphs with at least one edge \cite[Theorem~1.2]{Muta}. Ficarra, Moradi and Römer established componentwise linearity of symbolic powers for several classes of graphs with chordal complement \cite[Theorem~2.3]{FMR}. Ahmed and Namiq obtained $\reg(I^{(s)})=2s$ for every graph with chordal complement and at least one edge \cite[Theorem~1.1]{AN}.

If $G$ is \emph{gap-free} (no four vertices induce two disjoint edges) and has an edge, Theorem~\ref{thm:main} and Dupont García's ordinary-power formula \cite[Theorem~1]{DG} give
\[
 \reg(I^{(s)})=\reg(\overline{I^s})=\reg(I^s)
   =\max\{2s,\,\reg(I)+s-1\}\qquad(s\ge1).
\]
Minh and Vu established the formula for $s=2,3$ and conjectured it for both ordinary and symbolic powers for all $s\ge2$ \cite[Theorem~1.3 and Conjecture~A]{MVgap}.

Our proof bounds the regularity of $M_s=I^{(s)}/I^s$. Muta described the Krull dimension of this module in terms of induced odd cycles \cite[Corollary~4.10]{Muta}. We show that $M_s$, its submodule $\overline{I^s}/I^s$, and the corresponding quotient all have regularity at most $\reg(S/I^{(s)})$ (Theorem~\ref{thm:kernels}). The comparison of the three ideals then follows from two short exact sequences.

Our use of induced subgraphs and capacitated matchings in Section~\ref{sec:proof} is related to the odd-cycle bound of Kumar and Kumar \cite[Lemma~4.4]{KK} and to Dupont García's analysis of gap-free graphs \cite[Section~6]{DG}.

\section{Preliminaries}\label{sec:pre}

\subsection{Graphs, monomials and powers}\label{sub:powers}

For a finite simple graph $F$, let $R_F=\kk[x_v:v\in V(F)]$ and let
$I(F)=(x_ux_v:uv\in E(F))$ be its edge ideal. For $D\subseteq V(F)$,
write $F[D]$ for the subgraph induced by $D$, $F-D$ for the subgraph
induced by $V(F)\setminus D$, and $N_F[D]$ for $D$ together with its
neighbors. The open neighborhood of $v$ is
$N_F(v)=N_F[\{v\}]\setminus\{v\}$. Omitted graph subscripts refer to
$G$, and $\sqcup$ denotes disjoint union with no edges between the parts.

Write $\NN=\{0,1,\ldots\}$. For $a=(a_v)_{v\in V}\in\NN^V$, put
\[
 x^a=\prod_{v\in V}x_v^{a_v},\qquad
 P_a=\{v\in V:a_v>0\},\qquad
 a(D)=\sum_{v\in D}a_v \quad (D\subseteq V).
\]
Thus $x^a$ has multidegree $a$ and total degree $|a|=a(V)$. We write
$a|_D$ for the restriction of $a$ to $D$ and $\mathbf e_v$ for the
unit vector at $v$. For monomial ideals $K\subseteq L$, the residue
classes of the monomials in $L\setminus K$ form a $\kk$-basis of $L/K$.

Membership in a symbolic power of $I$ can be read from the vertex covers
of $G$. A \emph{vertex cover} is a set of vertices meeting every edge.
The minimal primes of $I$ are the ideals $\mathfrak p_C=(x_v:v\in C)$
indexed by the inclusion-minimal vertex covers $C$. For every integer
$s\ge1$,
\begin{equation}\label{eq:symbolic-pre}
 I^{(s)}
 =\bigcap_C \bigl(I^s S_{\mathfrak p_C}\cap S\bigr)
 =\bigcap_C \mathfrak p_C^{\,s},
\end{equation}
where both intersections run over the minimal covers and
$S_{\mathfrak p_C}$ denotes localization at $\mathfrak p_C$
\cite[proof of Proposition~1.5, p.~122]{HT}. In particular,
$x^a\in\mathfrak p_C^{\,s}$ precisely when $a(C)\ge s$.

The integral closure $\overline K$ of an ideal $K$ consists of the
elements $f$ satisfying $f^r+c_1f^{r-1}+\cdots+c_r=0$ for some
integer $r\ge1$ and coefficients $c_i\in K^i$. The integral closure
of a monomial ideal is again monomial \cite[Section~1.4.1]{HH}.

\subsection{Packings, covers and matchings}\label{sub:matching}

The exponent vector $a$ also has a packing interpretation. A product
of edge monomials divides $x^a$ precisely when each vertex $v$ occurs
among their endpoints at most $a_v$ times. Accordingly, an
\emph{edge packing} with capacities $a$ assigns an integer
$\lambda_e\ge0$ to each edge so that
$\sum_{e\ni v}\lambda_e\le a_v$ for every $v$. Its size is
$\sum_e\lambda_e$, and $\nu(a)$ denotes the maximum size. A fractional
packing allows nonnegative real edge weights. A packing \emph{saturates}
$v$ if $\sum_{e\ni v}\lambda_e=a_v$.

On the cover side, put $\tau(a)=\min_C a(C)$, where $C$ ranges over
all vertex covers of $G$. A fractional vertex cover is a point of
\[
 \mathcal C(G)
 =\{y\in[0,1]^V:y_u+y_v\ge1 \text{ for every }uv\in E\}.
\]
Its cost is $\sum_{v\in V}a_vy_v$, and $\sfq(a)$ denotes the minimum
cost. Subscripts on $\nu$, $\tau$ and $\sfq$ specify the graph when it
differs from $G$. All three optima are unchanged on restricting to
$G[P_a]$. Indeed, packings use only edges of $G[P_a]$, and covers
extend to $G$ at the same cost by assigning one outside $P_a$.

Fractional packing and covering are dual linear programs. The upper
bounds $y_v\le1$ preserve the dual optimum, since any coordinate
exceeding one can be lowered to one without increasing the cost.
Complementary slackness then relates the support of an optimal packing
to the coordinates of an optimal cover.

\begin{lem}[Packing--cover duality and complementary slackness
{\cite[Theorem~5.4 and (5.13)]{Sch}}]\label{lem:duality}
For $a\in\NN^V$, the maximum size of a fractional packing equals
$\sfq(a)$. For every optimal fractional packing $\lambda$ and every
minimum fractional cover $y$,
\begin{equation}\label{eq:slack}
 y_v>0 \ \Longrightarrow\ \sum_{e\ni v}\lambda_e=a_v,
 \qquad
 \lambda_{uv}>0 \ \Longrightarrow\ y_u+y_v=1.
\end{equation}
\end{lem}

The choice of an optimal cover can be refined further. Its coordinates
may be taken half-integral, and coordinates that are integral in
different optima can be made integral simultaneously.

\begin{lem}[Half-integrality and optimal covers
{\cite[Theorem~64.11 and (64.26)]{Sch}}]\label{lem:half-integrality}
Every extreme point of $\mathcal C(G)$ has coordinates in
$\{0,\tfrac12,1\}$ \cite[pp.~232--233]{NT}. For a fixed nonnegative
cost vector and two minimum fractional covers $y,y'$, there is a
minimum fractional cover $z$ that is integral at every vertex where
$y$ or $y'$ is integral.
\end{lem}

For a half-integral cover, a tight edge either joins two vertices of
value $1/2$ or joins a vertex of value zero to one of value one. The
edges of the latter type form a bipartite graph, where the packing
optimum is integral.

\begin{lem}[Bipartite packing integrality
{\cite[Theorem~21.2]{Sch}}]\label{lem:bipartite}
Let $F$ be bipartite with parts $L,R$ and integer capacities
$a\in\NN^{V(F)}$. There is an integral packing of size
$\tau_{\mathrm f,F}(a)$. In particular, if a fractional packing
saturates every vertex of $L$, then an integral packing does so as well.
\end{lem}

For general graphs, integer capacities can be represented by copies of
vertices. This reduces capacitated packing to ordinary matching. A
\emph{matching} is a set of vertex-disjoint edges. It is \emph{induced}
if its endpoints induce exactly those edges. For an unweighted graph
$F$, write $\nu(F)$ for its matching number and $\tau(F)$ for the
minimum cardinality of a vertex cover. A \emph{fractional perfect
matching} assigns nonnegative weights to the edges whose sum at each
vertex is one.

The \emph{replica graph} $G^a$ replaces each $v\in P_a$ by $a_v$
independent copies and each edge of $G[P_a]$ by all edges between the
corresponding copy classes. Projection sends each copy to its original
vertex.

\begin{lem}[Replica graph
{\cite[Proposition~2.3, p.~124]{HT}}]\label{lem:replica}
The matching and cover numbers of $G^a$ satisfy
\[
 \nu(G^a)=\nu(a),\qquad \tau(G^a)=\tau(a).
\]
The maximum fractional packing size in $G^a$ with unit capacities is
$\sfq(a)$. The weighted cover number satisfies $\tau(a)\le2\nu(a)$
\cite[proof of Theorem~3.3\textup{(iii)}, p.~127]{HT}.
\end{lem}

The number of vertices left uncovered by a maximum matching is its
\emph{deficiency}. The Tutte--Berge formula expresses this number in
terms of odd components.

\begin{lem}[Tutte--Berge formula
{\cite[Theorem~24.1]{Sch}}]\label{lem:tutte-berge}
For a finite graph $F$, write $\odd(F)$ for the number of its connected
components of odd order. Then
\begin{equation}\label{eq:tb-pre}
 |V(F)|-2\nu(F)
 =\max_{X\subseteq V(F)}\bigl(\odd(F-X)-|X|\bigr).
\end{equation}
\end{lem}

Returning to monomials, the three optima give parallel membership
criteria for the ordinary power, its integral closure, and the symbolic
power.

\begin{lem}[Monomial membership
{\cite[Propositions~1.1 and~1.5, pp.~121--122]{HT}}]\label{lem:membership}
For $a\in\NN^V$ and every integer $s\ge1$,
\[
\begin{aligned}
 x^a\in I^s            &\iff \nu(a)\ge s,\\
 x^a\in\overline{I^s}  &\iff \sfq(a)\ge s,\\
 x^a\in I^{(s)}        &\iff \tau(a)\ge s.
\end{aligned}
\]
\end{lem}

Since $\nu(a)\le\sfq(a)\le\tau(a)$, these criteria give
$I^s\subseteq\overline{I^s}\subseteq I^{(s)}$.

\subsection{Regularity and monomial quotients}\label{sub:reg}

To compare these ideals and their quotients, we fix our conventions
for regularity. Let $R=\kk[x_1,\ldots,x_n]$ with $\deg x_i=1$ for
every $i$. For a graded $R$-module $M$ and an integer $t$, the shift
$M(-t)$ is defined by $M(-t)_d=M_{d-t}$. If $M$ is nonzero and
finitely generated, write the $i$th module in its minimal graded free
resolution over $R$ as
\[
 F_i=\bigoplus_{j\in\mathbb Z}R(-j)^{\beta^R_{i,j}(M)}.
\]
Thus $\beta^R_{i,j}(M)$ counts the free generators of internal degree
$j$ in homological degree $i$. Following \cite[Section~4A]{Eis}, we
define the Castelnuovo--Mumford regularity of $M$ by
\[
 \reg_R(M)
 =\max\{j-i:i\ge0,\ j\in\mathbb Z,\ \beta^R_{i,j}(M)\ne0\}.
\]
The subscript specifies the ring of the resolution. Thus $\reg_S$
refers to the full vertex ring, whereas $\reg_{R_F}$ refers to the
polynomial ring on $V(F)$. We omit the subscript when the ring is clear
and set $\reg(0)=-\infty$. For a multigraded module we use the total
grading, so a multidegree $\eta\in\mathbb Z^n$ has internal degree
$|\eta|=\sum_{k=1}^n\eta_k$.

In particular, $\reg(M(-t))=\reg(M)+t$. For a nonzero proper
homogeneous ideal $K\subset R$, one has $\reg(K)=\reg(R/K)+1$.
For a graph $F$, put $h(F)=\reg(R_F/I(F))$. Thus $h(F)=0$ if $F$
is edgeless.

Regularity bounds pass from the factors of a finite filtration to the
whole module. Indeed, for a degree-preserving short exact sequence
$0\to U\to M\to W\to0$, the long exact sequence of graded
$\operatorname{Tor}$ gives
\begin{equation}\label{eq:reg-ses}
 \reg(M)\le\max\{\reg(U),\reg(W)\},\qquad
 \reg(W)\le\max\{\reg(M),\reg(U)-1\}.
\end{equation}
For a finite direct sum, regularity is the maximum of the regularities
of its summands.

We will also compare modules over different vertex rings. Extension
of scalars by polynomial variables preserves regularity. If $H$ is an
induced subgraph of $G$, a finitely generated graded $R_H$-module $M$
may instead be viewed as an $S$-module through the quotient map
$S\to R_H$ sending the variables outside $V(H)$ to zero. For this
action, $\reg_S(M)=\reg_{R_H}(M)$ by \cite[Corollary~4.6]{Eis}.

When one homogeneous generator is adjoined to an ideal, the resulting
cyclic quotient is described by a colon ideal. For a homogeneous ideal
$K\subseteq R$ and a homogeneous element $g$ of degree $e$, the ideal
$(K:g)=\{f\in R:fg\in K\}$ is the kernel of the map $f\mapsto fg+K$.
Hence
\begin{equation}\label{eq:cyclic-pre}
 (K+(g))/K\cong\bigl(R/(K:g)\bigr)(-e).
\end{equation}
For example, if $I=(xy,xz,yz)\subset\kk[x,y,z]$ is the edge ideal of
a triangle, then $I^{(2)}=I^2+(xyz)$ and $(I^2:xyz)=(x,y,z)$.
Thus $I^{(2)}/I^2\cong\kk(-3)$.

For monomial ideals, two operations on the defining ideal preserve
an upper bound on the regularity of its quotient: adjoining variables
and taking a colon by a monomial.

\begin{lem}[Monomial restrictions and colons
{\cite[Lemma~4.1\textup{(ii)} and Corollary~4.8]{Cav}}]\label{lem:restriction}
Let $K\subseteq R$ be a monomial ideal, $L$ an ideal generated by a
subset of the variables, and $p$ a monomial. Then
\[
 \reg(R/(K+L))\le\reg(R/K),\qquad
 \reg(R/(K:p))\le\reg(R/K).
\]
\end{lem}

For the first inequality, the cited corollary applies one variable at
a time after passing from a nonzero proper ideal to its quotient. If
$K=0$, then $R/L$ has regularity zero. If $K=R$, all the quotients
are zero.

\subsection{Regularity estimates for edge ideals}\label{sub:external}

For powers of an edge ideal, the colon ideals arising in
\eqref{eq:cyclic-pre} depend on the order in which the generators are
adjoined. Banerjee's ordering controls the contribution of each earlier
generator to the colon by the next one.

\begin{lem}[Banerjee's generator ordering
{\cite[Theorem~4.12, p.~310]{Ban}}]\label{lem:ordering}
Let $P$ be an edge ideal and $\ell\ge1$. Its power $P^\ell$ has an
ordering $g_1,\ldots,g_k$ of its minimal monomial generators with the
following property. Write $\rho_{ri}=g_r/\gcd(g_r,g_i)$, where
$\gcd$ denotes the greatest common monomial divisor. For each $r<i$,
either $\rho_{ri}\in(P^{\ell+1}:g_i)$ or some $\rho_{hi}$ with $h<i$
is a variable dividing $\rho_{ri}$.
\end{lem}

For a symbolic power, taking the colon by a product of edge monomials
admits a regularity bound in terms of the edge ideal itself.

\begin{lem}[Seyed Fakhari's colon bound
{\cite[Lemma~3.5 in the preprint]{Fak}}]\label{lem:symbolic-colon}
If $J$ is an edge ideal and $u$ is a product of $j\ge1$ of its edge
monomials, with repetitions allowed, then
\[
 \reg(J^{(j+1)}:u)\le\reg(J).
\]
\end{lem}

Lower bounds can be obtained by passing to an induced subgraph. When
that subgraph is a disjoint union, its components can be treated
separately, and the contribution of a matching is explicit.

\begin{lem}[Induced subgraphs and disjoint unions]\label{lem:graph-reg}
Let $s\ge1$ be an integer.
\begin{enumerate}
\item For every induced subgraph $F$ of $G$, the restriction inequality
\cite[Corollary~2.8]{MNP}, in quotient form, gives
\[
 \reg_{R_F}(R_F/I(F)^{(s)})\le\reg_S(S/I^{(s)}).
\]
If $F$ is edgeless, the left side is zero, and $S/I^{(s)}$ has a
minimal generator of degree zero.

\item Let $F=D\sqcup H$, with $D$ having at least one edge. If $H$
also has an edge, the term with exponents $(s,1)$ in the disjoint-sum
formula \cite[Theorems~5.6\textup{(ii)} and~5.11]{HNTT} gives
\[
 \reg_{R_F}(R_F/I(F)^{(s)})
 \ge \reg_{R_D}(R_D/I(D)^{(s)})+h(H).
\]
When $H$ is edgeless, equality holds by polynomial extension.

\item If $D$ consists of $m\ge1$ disjoint edges, its symbolic and
ordinary powers coincide \cite[Theorem~5.9]{SVV}, with
\[
 \reg(R_D/I(D)^{(s)})=2s+m-2
\]
by the formula for powers of forests \cite[Theorem~4.7]{BHT}.
\end{enumerate}
\end{lem}

\section{Proof of the main theorem}\label{sec:proof}

Fix $s\ge1$. The natural maps from $S/I^s$ and $S/\overline{I^s}$ onto $S/I^{(s)}$ reduce the problem to bounding the regularities of their kernels. We filter $I^{(s)}/I^s$ using minimum fractional vertex covers. In each summand of a layer, the exponents are fixed on the positive-weight vertices forced to take value $1/2$ and on their neighbors. The remaining induced graph $H$ contributes $h(H)$ to the upper bound. To obtain a lower bound for $\reg(S/I^{(s)})$ with the same term, we adjoin to $H$ an induced matching on the forced vertices. Its size accounts for the gap between fractional and integral packings. Since $\overline{I^s}/I^s$ is a term of the filtration, the comparison bounds both kernels.

\subsection{Forced half-integral vertices}\label{sec:half}

For $a\in\NN^V$, let $B(a)$ consist of the vertices in $P_a$ whose coordinate is $1/2$ in every minimum fractional cover of $G[P_a]$ with costs $a|_{P_a}$. Restriction to $P_a$ and extension by one on $V\setminus P_a$ show that the same set is obtained from optimal covers of $G$. Restricting to positive weights will matter when a capacity increases from zero to one. We first choose a single optimal cover that is integral off $B(a)$ and whose vertices of value one can all be saturated without using $B(a)$.

\begin{lem}\label{lem:canonical}
There is a half-integral optimal cover of $G[P_a]$ whose classes $O,B,Z$ of vertices with values $1,1/2,0$, respectively, satisfy the following properties.
\begin{enumerate}
\item $B=B(a)$, and the cover identically equal to $1/2$ on $B$ is the unique optimum on $G[B]$ with costs $a|_B$.
\item There is a fractional packing saturating every vertex of $B$ within $G[B]$, and an integral packing of size $a(O)$ on the edges between $O$ and $Z$, saturating $O$.
\end{enumerate}
\end{lem}
\begin{proof}
Work on $G[P_a]$ and put $B=B(a)$. For each $v\in P_a\setminus B$, some optimal extreme point is integral at $v$. Indeed, the optimal face is the convex hull of its extreme points, and these are half-integral by Lemma~\ref{lem:half-integrality}. Starting with any optimum and applying the second assertion of that lemma successively to the chosen extreme points gives an optimum integral on $P_a\setminus B$. On $B$ it is identically $1/2$ by definition. Let $O$ and $Z$ be its one- and zero-classes.

By Lemma~\ref{lem:duality}, an optimal fractional packing saturates $B\cup O$ and uses only edges within $B$ or between $O$ and $Z$. The latter edges form a bipartite graph, so Lemma~\ref{lem:bipartite} gives an integral packing saturating $O$. The packing within $B$ has size $a(B)/2$, equal to the cost of the all-half cover there. Any other minimum fractional cover of $G[B]$ could replace it in the global optimum, since every edge from $B$ to $P_a\setminus B$ ends in $O$. This would give a global optimum differing from $1/2$ on $B$, contrary to the definition of $B(a)$.
\end{proof}

Fix an optimal cover as in Lemma~\ref{lem:canonical}, with classes $O,B,Z$. Set $T=N[B]\setminus B$ and $H=G-N[B]$. Here the neighborhood is taken in the full graph even if some vertices in $T$ have zero weight. Write $b=a(B)$, $d=a(T)$, $\nu_B=\nu_{G[B]}(a|_B)$, and $\tau_B=\tau_{G[B]}(a|_B)$. All positive-weight vertices of $T$ lie in $O$. Put $j=a(O\setminus T)$. The packing between $O$ and $Z$ contributes equally to the three optima, so their differences are determined entirely by $G[B]$.

\begin{prop}\label{prop:split}
For the sets and weights above,
\begin{align}
 2\sfq(a)&=b+2(d+j),\label{eq:frac}\\
 \nu(a)&=\nu_B+d+j,\qquad \tau(a)=\tau_B+d+j.\label{eq:int}
\end{align}
Moreover, $j=\tau_H(a|_{V(H)})$. Counting multiplicities, the integral packing in Lemma~\ref{lem:canonical} uses $j$ edges in $H$ and $a_t$ edges from each $t\in T$ to $H$.
\end{prop}
\begin{proof}
The cost of the chosen cover gives \eqref{eq:frac}. Every edge of $G[P_a]$ outside $G[B]$ meets $O$, so at most $a(O)=d+j$ such edges occur in a packing. The integral $O$--$Z$ packing realizes this bound and uses no capacity on $B$. This proves the formula for $\nu$. An integral cover must pay at least $\tau_B$ on $B$ and at least $a(O)$ outside $B$ by the same packing. Taking a minimum cover on $B$, all of $O$, and all zero-weight vertices realizes equality.

All vertices of $Z$ belong to $H$. The $O$--$Z$ packing saturates $O$, and hence has the asserted multiplicities at $T$ and in $H$. Its $j$ internal edges give $\tau_H(a|_{V(H)})\ge j$. The set $O\setminus T$, together with the zero-weight vertices of $H$, covers $H$ at cost $j$.
\end{proof}

Define $q(a)=2\sfq(a)$, $b(a)=a(B(a))$, and $\pi(a)=(q(a),b(a))$, ordered lexicographically (compare $q$ first). Both entries are integers by \eqref{eq:frac}. Multiplication by $x_v$ increases the capacity at $v$ by one. To filter by $\pi$, we need to know how it changes under this operation and which exponents remain fixed when $\pi$ stays constant.

\begin{lem}\label{lem:persist}
For $v\in V$ let $a'=a+\mathbf e_v$. Then $\pi(a')\ge\pi(a)$. If $q(a')=q(a)$, then $B(a)\subseteq B(a')$ and $v\notin N[B(a)]$. If $\pi(a')=\pi(a)$, then also $B(a')=B(a)$ and $a'=a$ on $N[B(a)]$.
\end{lem}
\begin{proof}
For an optimal cover $y'$ of $G$ with costs $a'$,
\[
 \sfq(a')=\sum_{w\in V}a_w y'_w+y'_v\ge\sfq(a)+y'_v .
\]
Equality of the optima forces $y'_v=0$ and makes $y'$ optimal for $a$. Thus every new optimum is half on $B(a)$, and $v$ belongs neither to $B(a)$ nor to its neighborhood. The old forced vertices retain their weights, while each newly forced vertex has positive weight. Hence $b(a')\ge b(a)$, with equality only if the forced sets agree. The last assertion follows because only the coordinate at $v$ changed.
\end{proof}

\subsection{A filtration of the discrepancy module}\label{sec:filter}

Put $M_s=I^{(s)}/I^s$. Call $a\in\NN^V$ a \emph{discrepancy degree} if the class of $x^a$ in $M_s$ is nonzero. By Lemma~\ref{lem:membership}, this is equivalent to
\begin{equation}\label{eq:dis}
 \nu(a)<s\le\tau(a).
\end{equation}
Equation~\eqref{eq:int} then forces $B(a)\ne\varnothing$. Lemma~\ref{lem:replica} gives $\tau(a)\le2\nu(a)$, so
\[
 0\le b(a)\le q(a)\le 2\tau(a)\le4\nu(a)\le4(s-1).
\]
Thus $\pi$ takes finitely many values on the discrepancy degrees.

For a graph $H$, integers $j,d\ge0$, and subsets $U_1,\ldots,U_d$ of $V(H)$, possibly repeated, put $J=I(H)$ and $L_U=(x_v:v\in U)\subseteq R_H$. Define
\begin{equation}\label{eq:W}
 W_{j,d}(H;U_1,\ldots,U_d)
 =\frac{J^j\prod_{i=1}^d L_{U_i}+J^{(j+1)}}{J^{(j+1)}}.
\end{equation}
We take $J^0=R_H$ and the empty ideal product to be $R_H$.

Associate to each discrepancy degree the triple $\lambda=(B,\alpha,j)$, where $B=B(a)$, $\alpha=a|_{N[B]}$, and $j=\tau_H(a|_{V(H)})$ for $H=G-N[B]$. We consider only triples so obtained. Put $T=N[B]\setminus B$, $b=\alpha(B)$ and $d=\alpha(T)$. List each $t\in T$ exactly $\alpha_t$ times to obtain $t_1,\ldots,t_d$. Set $U_i=N_G(t_i)\cap V(H)$ and write $W_\lambda$ for \eqref{eq:W} with these sets. Their order does not affect the product. By \eqref{eq:frac}, $\pi(a)$ is determined by $\lambda$ and equals $\pi_\lambda=(b+2(d+j),b)$. The product in \eqref{eq:W} imposes the capacities on $H$ simultaneously on the $j$ internal edges and the $d$ edges from $T$ to $H$.

\begin{lem}\label{lem:fiber}
Fix such a triple $\lambda$ and $c\in\NN^{V(H)}$. Let $a'$ agree with $\alpha$ on $N[B]$ and with $c$ on $V(H)$. Then $a'$ is a discrepancy degree with triple $\lambda$ if and only if
\[
   x^c\in J^j\prod_{i=1}^{d}L_{U_i}\ \setminus\ J^{(j+1)} .
\]
\end{lem}
\begin{proof}
If $a'$ has triple $\lambda$, the packing in Proposition~\ref{prop:split} uses $j$ edges of $H$ and one edge from each listed $t_i$ to a vertex of $U_i$. Their endpoints in $H$, counted with multiplicity, respect the capacities $c$. This proves the required divisibility. Since $\tau_H(c)=j$, Lemma~\ref{lem:membership} gives $x^c\notin J^{(j+1)}$.

Conversely, divisibility gives $j$ edges of $H$ and a choice of one endpoint in each $U_i$, all within the capacities $c$. Joining the latter vertices to the corresponding $t_i$ gives an integral packing of $d+j$ edges outside $B$. Its $j$ internal edges imply $\tau_H(c)\ge j$, while $x^c\notin J^{(j+1)}$ gives the reverse inequality. Extend the indicator of a minimum-weight vertex cover of $H$ by assigning one on $T$ and one half on $B$. No edge joins $B$ to $H$, so this is a fractional cover of cost $b/2+d+j$.

Because $\lambda$ occurs for some discrepancy degree and $\alpha$ is fixed, the packing on $G[B]$ in Lemma~\ref{lem:canonical} still has size $b/2$ and saturates $B$. For any fractional cover, the internal packing forces a cost of at least $b/2$ on $B$, and the exterior packing forces a cost of at least $d+j$ on $V\setminus B$. The extended cover attains their sum. Every optimum therefore attains the internal bound, so its restriction to $B$ is the unique all-half optimum of Lemma~\ref{lem:canonical}. The extended cover is integral at every vertex outside $B$. Hence $B(a')=B$. Its residual cover number is $j$, and \eqref{eq:int} now gives the same $\nu$ and $\tau$ as for the original degree with triple $\lambda$. Thus \eqref{eq:dis} holds for $a'$.
\end{proof}

Lemma~\ref{lem:fiber} identifies the monomials with a fixed triple. By Lemma~\ref{lem:persist}, multiplication that survives at a fixed value of $\pi$ preserves that triple. This gives a direct sum decomposition of each filtration layer.

\begin{prop}\label{prop:filtration}
For an attained pair $\gamma$, let $F_{\ge\gamma}$ and $F_{>\gamma}$ be the $\kk$-spans in $M_s$ of monomial classes with $\pi(a)\ge\gamma$ and $\pi(a)>\gamma$, respectively. These spaces define a finite filtration by graded $S$-submodules, with factors
\begin{equation}\label{eq:layer}
 F_{\ge\gamma}/F_{>\gamma}\ \cong\
 \bigoplus_{\lambda:\,\pi_\lambda=\gamma} W_\lambda(-|\alpha|).
\end{equation}
In each summand $|\alpha|=b+d$, and the variables indexed by $N[B]$ act by zero. The submodule $\overline{I^s}/I^s$ is spanned by the classes of $x^a$ in $M_s$ at discrepancy degrees satisfying $q(a)\ge2s$, and $I^{(s)}/\overline{I^s}$ inherits the factors with $q(a)<2s$.
\end{prop}
\begin{proof}
Lemma~\ref{lem:persist} shows that multiplication by a variable preserves the indicated spans. For fixed $\gamma$ there are finitely many triples: $B\subseteq V$, $|\alpha|=b+d\le q(a)$, and $j$ is determined by $q(a)$, $b$ and $d$. A multiplication that survives in a factor preserves $B$ and $\alpha$ by the same lemma, and hence preserves $j$ by \eqref{eq:frac}.

For each triple, send the class of $x^c$ to the class of $x^\alpha x^c$, extending both exponent vectors by zero outside their respective vertex sets. Lemma~\ref{lem:fiber} identifies the monomial bases and gives the shift $|\alpha|$. Multiplication by any $x_v$ with $v\in V(H)$ agrees on the two sides: the product remains in the ideal product in \eqref{eq:W}, and it survives in the same factor exactly when it stays outside $J^{(j+1)}$. A variable indexed by $N[B]$ kills the image in the factor by Lemma~\ref{lem:persist}. This proves \eqref{eq:layer} as an isomorphism of singly graded $S$-modules. Finally, Lemma~\ref{lem:membership} identifies the integral-closure degrees with $q(a)\ge2s$. These are an upper interval of the filtration.
\end{proof}

\subsection{Regularity on the residual graph}\label{sec:local}

Proposition~\ref{prop:filtration} reduces the upper bound to a regularity estimate for $W_{j,d}(H;U_1,\ldots,U_d)$. The product ideal in the numerator of \eqref{eq:W}, when nonzero, is generated in degree $2j+d$. Realizing these generators as edge products in a larger graph lets us apply Banerjee's ordering to bound the regularity of their images modulo $J^{(j+1)}$.

\begin{thm}\label{thm:local}
For every finite simple graph $H$, all integers $j,d\ge0$, and all subsets $U_1,\ldots,U_d$ of $V(H)$,
\[
  \reg\bigl(W_{j,d}(H;U_1,\ldots,U_d)\bigr)\le h(H)+2j+d .
\]
\end{thm}
\begin{proof}
Set $R=R_H$, $J=I(H)$, $Q=J^{(j+1)}$ and $\ell=j+d$. The case $\ell=0$ gives $W_{0,0}=R/J$. Assume $\ell>0$ and the module is nonzero.

\smallskip\noindent\emph{Realization in one multidegree.}
Adjoin independent vertices $z_1,\ldots,z_d$ to $H$, with $z_i$ adjacent precisely to $U_i$, and denote the edge ideal of this graph by $P\subset R'=R[z_1,\ldots,z_d]$. We use $z_i$ also for the corresponding variable. Give $R'$ the auxiliary grading by the exponents of the $z_i$, as well as the total degree in the old variables. Write $Z=z_1\cdots z_d$ and $\mathbf1=(1,\ldots,1)\in\NN^d$. For $d=0$ these conventions give $Z=1$ and a trivial auxiliary grading. Division by $Z$ in the component of degree $\mathbf1$ identifies
\[
 \bigl((P^\ell+QR')/QR'\bigr)_{\mathbf1}
 \ \cong\ W_{j,d}(H;U_1,\ldots,U_d)
\]
as graded $R$-modules, with degrees measured in the old variables. Indeed, a product of $\ell$ edges dividing $Zx^c$ uses each new vertex at most once. If it uses fewer than $d$, at least $j+1$ edges lie in $H$, so $x^c\in J^{j+1}\subseteq Q$. If it uses all $d$, the old variables contribute a product of $j$ edges of $H$ and one variable from each $U_i$. Conversely, every such product occurs.

\smallskip\noindent\emph{The cyclic factors.}
Order the minimal generators $g_1,\ldots,g_k$ of $P^\ell$ as in Lemma~\ref{lem:ordering}. Set $E_0=QR'$ and $E_i=E_{i-1}+(g_i)$. Taking the component of auxiliary degree $\mathbf1$ preserves exact sequences. A generator with a repeated new variable contributes nothing to that component. A generator using fewer than $d$ new variables, each at most once, belongs to $QR'$. Hence $(E_i/E_{i-1})_{\mathbf1}\ne0$ requires $g_i=Zg$, with $g\in R$ of degree $2j+d$. By \eqref{eq:cyclic-pre} this component is
\begin{equation}\label{eq:cyclic}
 (E_i/E_{i-1})_{\mathbf1}\cong(R/D_i)(-(2j+d)),
 \qquad D_i=(E_{i-1}:g_i)\cap R.
\end{equation}
Indeed, all multipliers of $g_i$ that preserve its auxiliary degree lie in $R$.

Put $L_i=(x_v:v\in V(H),\ x_v\in D_i)$. We claim
\begin{equation}\label{eq:colon}
 D_i=(Q:g)+L_i .
\end{equation}
The inclusion from right to left is immediate. For the reverse inclusion, let $w\in D_i\setminus(Q:g)$ be a monomial. Some $g_r$ with $r<i$ divides $wg_i$, so $\rho_{ri}=g_r/\gcd(g_r,g_i)$ divides $w$. If $\rho_{ri}\in(P^{\ell+1}:g_i)$, then $\ell+1$ edges can be packed into $wg_i=Zwg$. At most $d$ of them use a new vertex, whence $wg\in J^{j+1}\subseteq Q$, a contradiction. Banerjee's other alternative gives a variable $\rho_{hi}$ dividing $w$ for some $h<i$. It lies in $R$ and its product with $g_i$ is divisible by $g_h$, so it belongs to $L_i$. This proves \eqref{eq:colon}.

\smallskip\noindent\emph{The regularity bound.}
Write $g=up$, where $u$ is a product of $j$ edge monomials of $H$ and $p$ is a product of one variable from each $U_i$. Since the factor \eqref{eq:cyclic} is nonzero, $g\notin Q$ and $(Q:u)$ is proper. If $j>0$, then $H$ has an edge and $(Q:u)$ contains $Q\ne0$. Lemma~\ref{lem:symbolic-colon} in quotient form gives $\reg(R/(Q:u))\le h(H)$. If $j=0$, this bound is an equality because $u=1$ and $Q=J$. Applying Lemma~\ref{lem:restriction} first to the colon by $p$ and then to the variables in $L_i$, we obtain from \eqref{eq:colon}
\[
 \reg(R/D_i)\le\reg(R/(Q:g))\le\reg(R/(Q:u))\le h(H).
\]
After the shift in \eqref{eq:cyclic}, each factor has regularity at most $h(H)+2j+d$. The filtration bound \eqref{eq:reg-ses} completes the proof.
\end{proof}

\subsection{An induced matching and the regularity comparison}\label{sec:last}

For a discrepancy degree $a$, put $B=B(a)$ and $H=G-N[B]$. The complementary lower bound calls for an induced matching in $G[B]$ that can be adjoined to $H$. The number of vertices left uncovered by a maximum matching in the replica graph of $G[B]$ determines how many edges we need.

\begin{lem}\label{lem:deficiency}
Let $a$ be a discrepancy degree, $B=B(a)$, and $\delta=a(B)-2\nu_B$. Then $G[B]$ contains an induced matching of size $\max\{1,\delta\}$.
\end{lem}
\begin{proof}
Equations \eqref{eq:dis} and \eqref{eq:int} give $B\ne\varnothing$. Let $F$ be the replica graph of $G[B]$ with capacities $a|_B$. It has $a(B)$ vertices and matching number $\nu_B$ by Lemma~\ref{lem:replica}. The fractional packing of Lemma~\ref{lem:canonical} saturates $B$, so its size is $a(B)/2$. By Lemmas~\ref{lem:duality} and~\ref{lem:replica}, $F$ admits a fractional matching of this size. It is perfect because the sum of all incident weights, counted over the vertices of $F$, is $a(B)=|V(F)|$.

If $\delta=0$, then $F$ has an edge, whose projection suffices. Suppose $\delta>0$. By Lemma~\ref{lem:tutte-berge}, choose $X\subseteq V(F)$ with $\odd(F-X)=|X|+\delta$. Each singleton component of $F-X$ sends one unit of fractional matching weight to $X$. The vertices of $X$ receive at most $|X|$ units in total, so $F-X$ has at most $|X|$ singleton components. At least $\delta$ of its odd components therefore contain an edge. Choose one edge in each of $\delta$ such components. Their endpoints induce precisely the chosen edges in $F$.

The projected edges also form an induced matching in $G[B]$. If two shared a vertex $v$, the copy of $v$ on either chosen edge would be adjacent to the opposite endpoint of the other. Likewise, any edge joining endpoints of distinct projected edges would lift to an edge between the chosen components. Both conclusions contradict the definition of the components of $F-X$.
\end{proof}

The matching from Lemma~\ref{lem:deficiency} and the residual graph $H$ induce a disjoint union in $G$. Lemma~\ref{lem:graph-reg} therefore gives a lower bound for $\reg(S/I^{(s)})$ with the same term $h(H)$ as the upper bound from Theorem~\ref{thm:local}. Comparing the remaining terms gives the following estimate for all three discrepancy modules.

\begin{thm}\label{thm:kernels}
Let $J_s=\overline{I^s}$ and $A_s=S/I^{(s)}$. Each of the modules
\[
 I^{(s)}/I^s,\qquad I^{(s)}/J_s,\qquad J_s/I^s
\]
has regularity at most $\reg(A_s)$.
\end{thm}
\begin{proof}
Fix a summand of \eqref{eq:layer} and a discrepancy degree $a$ with triple $\lambda=(B,\alpha,j)$. Use $H,b,d$ as in Proposition~\ref{prop:split}. Since the variables of $N[B]$ act by zero, Theorem~\ref{thm:local} and \eqref{eq:frac} give
\begin{equation}\label{eq:upper}
 \begin{aligned}
 \reg_S\bigl(W_\lambda(-|\alpha|)\bigr)
 &= |\alpha|+\reg_{R_H}(W_\lambda) \\
 &\le (b+d)+h(H)+2j+d=q(a)+h(H).
 \end{aligned}
\end{equation}
Let $\delta=b-2\nu_B$ and choose the induced matching of size $m=\max\{1,\delta\}$ in $G[B]$ given by Lemma~\ref{lem:deficiency}. No vertex of this matching is adjacent to a vertex of $H=G-N[B]$. Writing $D$ for the graph induced by the matching endpoints, we thus have an induced subgraph $F=D\sqcup H$ of $G$. Lemma~\ref{lem:graph-reg} gives
\[
 \reg_S(A_s)\ge\reg_{R_F}(R_F/I(F)^{(s)})
 \ge\reg_{R_D}(R_D/I(D)^{(s)})+h(H)
 =2s-2+m+h(H).
\] Subtracting twice the first equality of \eqref{eq:int} from \eqref{eq:frac} gives $q(a)-2\nu(a)=\delta$. Since $\nu(a)\le s-1$,
\begin{equation}\label{eq:lower}
 q(a)+h(H)\le 2s-2+m+h(H)\le\reg(A_s) .
\end{equation}
Together, \eqref{eq:upper} and \eqref{eq:lower} bound every summand by $\reg(A_s)$. Taking finite direct sums and successive extensions gives the assertion for $M_s$. By Proposition~\ref{prop:filtration}, the same argument applies to the subfiltration with $q\ge2s$ and the quotient filtration with $q<2s$, including the cases in which one is empty.
\end{proof}

\begin{proof}[Proof of Theorem~\ref{thm:main}]
For either $D=I^s$ or $D=\overline{I^s}$, Theorem~\ref{thm:kernels} bounds the regularity of the kernel in the exact sequence
\[
 0\longrightarrow I^{(s)}/D\longrightarrow S/D\longrightarrow A_s\longrightarrow0 .
\]
Put $r=\reg(A_s)$. The kernel bound and \eqref{eq:reg-ses} give
\[
 \reg(S/D)\le r,\qquad
 r\le\max\{\reg(S/D),r-1\}.
\]
As $A_s$ is a nonzero quotient of $S$, $r$ is finite. Thus $\reg(S/D)=r$ for both choices of $D$. When $G$ has an edge, passing from each quotient to its defining ideal adds one to regularity. When $G$ is edgeless, all three ideals are zero.
\end{proof}

\begin{ex}\label{ex:four}
Let $G$ be the disjoint union of two triangles and two single edges. Both inclusions $I^4\subsetneq\overline{I^4}\subsetneq I^{(4)}$ are strict. In each of the following assignments, all unlisted vertices have weight zero. Unit weights on one triangle and both single edges give $(\nu(a),\sfq(a),\tau(a))=(3,7/2,4)$. Unit weights on both triangles and one single edge give $(3,4,5)$. In the latter case $B(a)$ consists of the six triangle vertices, $\nu_B=2$, and $\delta=2$. One edge from each triangle gives the induced matching in \eqref{eq:lower}.
\end{ex}

\section*{Declaration on the use of generative AI}
The authors used ChatGPT-6.2 Astra to assist in discussing proof strategies, generating proofs, checking proofs, and polishing the paper. The authors also used ChatGPT-6.2 Sol Codex to perform Lean 4 formal verification of the paper. The main theorem and its dependencies have been re-checked by the Lean kernel, and an axiom audit found no proof placeholders or custom axioms.
The authors take full responsibility for the correctness
of the paper.

\bigskip

{\it E-mail address:} chenhao@suda.edu.cn

\medskip

{\it E-mail address:} 20254207033@stu.suda.edu.cn

\end{document}